\documentclass{jag-l}
\input xy
\xyoption{all}

\newcommand{\N}{{\mathbb N}}
\newcommand{\Z}{{\mathbb Z}}
\newcommand{\C}{{\mathbb C}}
\newcommand{\Q}{{\mathbb Q}}
\renewcommand{\P}{{\mathbb P}}
\newcommand{\p}{{\mathbb P}}
\def\g{{\mathfrak g}}
\def\W{\mathcal W}
\def\h{\mathfrak{h}}
\def\iso{\cong}
\DeclareMathOperator{\sgn}{sgn}
\newcommand{\gbinom}[2]{ {{#1}\brace{#2}} }

\newtheorem{theorem}{Theorem}[section]
\newtheorem{lemma}[theorem]{Lemma}
\newtheorem{proposition}[theorem]{Proposition}
\newtheorem{corollary}[theorem]{Corollary}

\theoremstyle{definition}

\newtheorem{example}[theorem]{Example}

\theoremstyle{remark}
\newtheorem{remark}[theorem]{Remark}

\numberwithin{equation}{section}

\begin{document}

\title[Degree of the discriminant]{The degree of the discriminant of irreducible representations}

\author{L. M. Feh\'er}
\address{Department of Analysis, Eotvos University Budapest, Hungary}
\email{lfeher@renyi.hu}
\author{A. N\'emethi}
\address{Renyi Mathematical Institute Budapest, Hungary;  and Ohio State University, Columbus OH, USA}
\email{nemethi@renyi.hu and nemethi@math.ohio-state.edu}
\author{R. Rim\'anyi}
\address{Department of Mathematics, University of North Carolina at Chapel Hill, USA}
\email{rimanyi@email.unc.edu}

\thanks{\noindent Supported by NSF grant DMS-0304759, and OTKA 42769/46878 (2nd author),
        NSF grant DMS-0405723 (3rd author) and OTKA T046365MAT (1st and 3rd author)
        }






\begin{abstract} We present a formula for the degree of the discriminant of
irreducible representations of a Lie group, in terms of the roots of
the group and the highest weight of the representation. The proof
uses equivariant cohomology techniques, namely, the theory of Thom
polynomials, and a new method for their computation. We study the
combinatorics of our formulas in various special cases.
\end{abstract}

\maketitle

\section{Introduction}

Let $G$ be a complex connected reductive algebraic group, and let
$\rho:G\to GL(V)$ be an irreducible algebraic representation. Then
$\rho$ induces an action of $G$ on the projective space $\P(V)$.
This action has a single closed orbit, the orbit of the weight
vector of the highest weight $\lambda$. E.g., for $GL(n)$ acting on
$\Lambda^k\C^n$, we get the Grassmannian $Gr_k(\C^n)$. The dual $\P
D_\lambda$ of this orbit (or the affine cone $D_\lambda$ over it) is
called the {\em discriminant} of $\rho$ since it generalizes the
classical discriminant. The goal of the present paper is to give a
formula for the degree of the discriminant in terms of the highest
weight of the representation $\rho$ and the roots of $G$ (Theorems
\ref{short} and \ref{long}).

The classical approach to find the degree of dual varieties is due
to Kleiman \cite{kleiman} and Katz \cite{katz}. Their method,
however, does not produce a formula in the general setting. Special
cases were worked out by Holme \cite{holme}, Lascoux \cite{lascoux},
Boole, Tevelev \cite{tevelevcikk}, Gelfand-Kapranov-Zelevinsky
\cite[Ch.13,14]{gkz}, see a summary in \cite[Ch.7]{tevelev}.
De~Concini and Weyman \cite{weyman} showed that, if $G$ is fixed,
then for regular highest weights the formula for the degree of the
discriminant is a polynomial with positive coefficients, and they
calculated the constant term of this polynomial. A corollary of our
result is an explicit form for this polynomial (Cor.~\ref{fehpol})
with the additional  fact that the same polynomial calculates the
corresponding degrees for non-regular highest weights as well
(modulo an explicit factor).

In special cases our formula can be expressed in terms of some basic
concepts in the combinatorics of polynomials
\cite{macdonald,lascouxbook}, such as the Jacobi symmetrizer,
divided difference operators, or the scalar product on the space of
polynomials (Section \ref{combinat}).

For the group $GL(n)$ we further simplify the formula in many
special cases in Section \ref{andras}.

The authors thank M. Kazarian, A. Knutson, S. Kumar, and A. Szenes
for helpful discussions. The original proof of Theorem \ref{main}
was more complicated. We got the idea of a simpler proof from a
lecture of A. Szenes.  We thank the referees for valuable
suggestions.

\section{Degree and Thom polynomials} \label{sec:deg&thom}

In this paper we will use cohomology with rational coefficients. The
Lie groups we consider are complex connected reductive algebraic
groups. All varieties are over the complex numbers $\C$ and
$GL(n)=GL(n;\C)$ denotes the general linear group of $\C^n$.

\subsection{Degree and cohomology} \label{sec:deg&coh}
Suppose that $Y$ is a smooth complex algebraic variety and $X\subset
Y$  is a closed subvariety of complex codimension $d$. Then $X$
represents a cohomology class $[X]$ in the cohomology group
$H^{2d}(Y)$. This class is called the {\em Poincar\'e dual} of $X$.
The existence of this class and its basic properties are explained
e.g. in \cite{fulton:young}. If $Y$ is the projective space $\P^n$
then $H^*(Y)\iso \Q[x]/(x^{n+1})$ and $[X]=\deg(X)x^d$, where $x$ is
the class represented by a hyperplane. By definition, the cone
$CX\subset \C^{n+1}$ of $X$ has the same degree.

\subsection{Degree and equivariant cohomology} \label{sec:deg&equi} We
would also like to express the degree in terms of  {\em equivariant
cohomology}. Let $G$ be a complex connected reductive algebraic
group (though some definitions and claims hold for more general
groups as well). Let $G$ act on a topological space $Y$. Then the
equivariant cohomology ring $H^*_G(Y)$ is defined as the ordinary
cohomology of the  Borel construction $EG\times_G Y$. Here $EG$
denotes the universal principal $G$-bundle over the classifying
space $BG$. If $Y$ is a smooth complex algebraic variety, $G$ is a
complex Lie group acting on $Y$ and $X\subset Y$ is a $G$-invariant
subvariety of complex codimension $d$ then  $X$ represents an {\em
equivariant} cohomology class (sometimes called equivariant
Poincar\'e dual or, if $Y$ is contractible, Thom polynomial) $[X]\in
H^{2d}_G(Y)$. This class has the following universal property.

Let $P\to M$ be an algebraic principal $G$-bundle. Then $P$ is
classified by a map $k: M\to BG$, in other words we have the diagram
\[ \xymatrix{P\ar[r]^{\tilde k}\ar[d]&EG\ar[d]\\M\ar[r]^{k}&BG.} \]
We also have an associated bundle $P\times_G Y\to M$ and an induced
map $\hat k:P\times_G Y\to EG\times_G Y$. The universal property of
$[X]\in H^{2d}_G(Y)$ is that for {\em every} $P\to M$ the ordinary
cohomology class $[P\times_G X]$ represented by the subvariety
$P\times_G X$ of $P\times_G Y$ satisfies
\begin{equation} \label{epd-def} \hat k^*[X]=[P\times_G X]. \end{equation}
It is easy to see that (\ref{epd-def}) characterizes the equivariant
cohomology class $[X]\in H^*(BG\times_G Y)=H_G^*(Y)$ uniquely. The
existence of such a class is not so obvious, for details see
\cite{kaza} or \cite{fr}.

Let  $Y=V$ be a complex vector space and  let $G=GL(1)$  act as
scalars. If $X\subset Y$ is a $G$-invariant  subvariety  of complex
codimension $d$ then $X$ is the cone of a projective variety $\P
X\subset \P V$. Then $H^{*}_G(V)\iso \Q[x]$ and $[X]=\deg(\P X)x^d$,
where $x$ is the class of a hyperplane. For more general complex Lie
groups we can expect to extract the degree if $G$ `contains the
scalars'. As we will see, even this condition is not necessary:  we
can simply replace $G$ with $G\times GL(1)$.

Suppose that \( \rho:G\to GL(V) \) `contains the scalars' i.e. there
is a homomorphism $h:GL(1)\to G$, such that $\rho\circ h$ is the
scalar representation on $V$: for $v\in V$ and $z\in GL(1)$ we have
$\rho\circ h(z)v=zv$. By the basic properties of the equivariant
Poincar\'e dual we have a `change of action' formula:

\begin{equation} [X]_{\rho\circ h}=h^*[X]_{\rho},\label{eq:change} \end{equation}
where
$$h^*:H^*_G(V)\iso H^*(BG) \to H^*_{GL(1)}(V)\iso H^*(BGL(1))$$
is induced by the map $Bh:BGL(1)\to BG$ classifying the principal
$G$-bundle $EGL(1)\times_h G$. The $\rho$ and $\rho\circ h$ in the
lower index indicates whether we take the $G$- or
$GL(1)$-equivariant Poincar\'e dual of $X$. The  `change of action'
formula implies that
    $$\deg(X)x^d=h^*[X]_{\rho}.$$
It frequently happens that we can only find a homomorphism
$h:GL(1)\to G$, such that $\rho\circ h(z)v=z^kv$ for some non zero
integer $k$. Then, by the same way we obtain
     $$k^d\deg(X)x^d=[X]_{\rho\circ h}=h^*[X]_{\rho}.$$
The calculation of $h^*$ is fairly simple. Suppose that
$m:GL(1)^r\to G$ is a (parametrized) maximal complex torus of $G$.
Then by Borel's theorem \cite[\S 27]{borel-th} $H^*(BG)$ is
naturally isomorphic to the Weyl-invariant subring of $H^*(BT)$ ($T$
is the image of $m$, a complex maximal torus of $G$), and $H^*(BT)$
can be identified with the symmetric algebra of the character group
of $T$. Hence $H^*(BG)=\Q[\alpha_1,\dotsc,\alpha_r]^W$, where the
$\alpha_i$'s generate the weight lattice of $G$ (one can identify
$\alpha_i$ with the $i$\textsuperscript{th} projection
$\pi_i:GL(1)^r\to GL(1)$) and $W$ is the Weyl group of $G$.  Then
$[X]=p(\alpha_1,\dotsc,\alpha_r)$ for some homogeneous polynomial
$p$  of degree $d$. We can assume that the homomorphism $h:GL(1)\to
G$ factors through $m$, i.e. $h=m\circ \phi$, where
$\phi(z)=(z^{k_1},\dotsc,z^{k_r})$ for some integers
$k_1,\dotsc,k_r$. Applying  the `change of action' formula
(\ref{eq:change}) once more leads~to

\begin{proposition} \label{deg} For the polynomial $p$, and integers $k,k_1,\dotsc,k_r$ defined
above $$\deg(X)=p(k_1/k,\dotsc,k_r/k).$$  \end{proposition}

This innocent-looking statement provides a uniform approach for
calculating the degree of degeneracy loci whenever a Chern--class
formula is known. A similar but more involved argument (see
\cite{forms}) calculates the equivariant Poincar\'e dual of $\P X$
in $\P V$.


\section{Calculation of the equivariant Poincar\'e dual} \label{sec:calc}

\label{sec:push} Now we reduce the problem of computing an
equivariant Poincar\'e dual to computing an integral. We will need
$G$-equivariant characteristic classes: the equivariant Chern
classes $c_i(E)\in H^{2i}_G(M)$ and the equivariant Euler class
$e(E)\in H^{2n}_G(M)$ of a $G$-equivariant complex vector bundle
$E\to M$ of rank $n$ are defined via the Borel construction. We have
$e(E)=c_{n}(E)$.

As in the case of ordinary cohomology, pushforward can be defined in
the equivariant setting. An introduction can be found in
\cite{atiyah-bott}. Its properties are similar; for example if
$\phi:\tilde X\to Y$ is a $G$-equivariant resolution of the
$d$-codimensional invariant subvariety $X\subset Y$, then
$[X]=\phi_*1\in H^{2d}_G(Y)$.

Let $V$ be a vector space and $X\subset V$ a subvariety. Suppose
that $M$ is a compact manifold and let $\tilde X \to M$ be a
sub-vector bundle $S$ of $M\times V \to M$. Let $\pi: M\times V \to
V$ be the projection, $i:\tilde X \subset M\times V$ the embedding,
and let $\phi=\pi \circ i$, as in the diagram
$$\xymatrix{ \tilde{X} \ar@/^2pc/[rrr]^{\phi} \ar@{^{(}->}@<-3pt>[rr]^{i}  \ar[rd] &  & M\times
V \ar[ld]
\ar[r]^{\pi} & V \\
 & M &  & X. \ar@{^{(}->}[u] } $$


\begin{proposition}\label{31} Suppose that $G$ acts on all spaces in the diagram above, and that all maps are $G$-invariant.
Let $Q$ be the quotient bundle $V\ominus S$ over $M$. We have
   \begin{equation} \pi_*i_*1=\int_{M}e(Q)\qquad \in H_G^*(V)=H^*(BG).\label{eq:intnonres} \end{equation}
If $\phi$ is a resolution of $X\subset V$ then
  \begin{equation}[X]=\int_{M}e(Q)\qquad \in H_G^*(V)=H^*(BG).\label{eq:int} \end{equation}
\end{proposition}

\begin{proof}
Observe that $i_*1$ is the Euler class of $Q=V\ominus S$ (we
identify the cohomology of $M\times V$ with the cohomology of $M$);
and $\pi_*$ is integration along the fiber of $\pi$.
\end{proof}

\begin{remark}\label{intzero} If $\phi(\tilde X)$ has smaller dimension than $\tilde X$ then $\pi_*i_*1=\int_{M}e(Q)$ is zero since this cohomology class is supported on $\phi(\tilde X)$ and its codimension is bigger than the rank of $\pi_*i_*1$.
\end{remark}

Below we will calculate the integral (\ref{eq:int}) using the
Berline-Vergne-Atiyah-Bott integral formula \cite{BV},
\cite{atiyah-bott}, that we recall now.

\begin{proposition}[Berline-Vergne-Atiyah-Bott] \label{ab} Suppose that $M$ is
  a compact oriented manifold,
$T$ is a torus acting smoothly on $M$ and $C(M)$ is the set of
components of the fixed point manifold. Then for any cohomology
class $\alpha\in H_T^*(M)$
\begin{equation}\label{abformula}\int_M\alpha=\sum_{F\in C(M)}
\int_F \frac{i^*_F\alpha}{e(\nu_{F})}.\end{equation}
 Here $i_F:F\to
M$ is the inclusion, $e(\nu_{F})$ is the $T$-equivariant Euler class
of the normal bundle $\nu_{F}$ of $F\subset M$. The right side is
considered in the fraction field of the polynomial ring of
$H^*_T($point$)=H^*(BT)$ (see more on details in
\cite{atiyah-bott}): part of the statement is that the denominators
cancel when the sum is simplified.
\end{proposition}

\begin{remark}
  Here  $T$ denotes a {\em real} torus $T=U(1)^r$. Since $BU(1)^r$ is homotopy equivalent to $BGL(1)^r$ we can always restrict the complex torus action to the real one without losing any cohomological information.
\end{remark}

Returning to the situation of Proposition~\ref{31}, let us further
assume that the group $G$ is a torus, and that the fixed point set
$F(M)$ on the compact manifold $M$ is finite. For a fixed point
$f\in F(M)$ the Euler class $e(V)$ restricted to $f$ can be computed
as the product of the Euler classes of the restrictions of $S$ and
$Q$ to $f$, that is $e(V_f)=e(S_f)e(Q_f)$. Also, $\nu_f$ is just the
tangent space $T_f M$. Hence, from Proposition~\ref{31} and \ref{ab}
we obtain that
\begin{equation}\label{we}
[X]=e(V)\sum_{f\in F(M)}\left(e(S_f)\cdot e(T_f M)\right)^{-1}.
\end{equation}
This formula is a special case of \cite[Prop. 3.2]{berczi-szenes}.
These type of formulas have their origin in \cite{Vergne}.

\section{The equivariant cohomology class of the dual of smooth varieties}

Let the torus $T$ act on the complex vector space $V$ and let
$X\subset V$ be a $T$-invariant cone. Assume moreover that $\P
X\subset \P V$ is smooth, and that the projective dual $\P\check X$
of $\P X$ is a hypersurface. As a general reference we refer the
reader to the book of Tevelev \cite{tevelev} on Projectively Dual
Varieties. For a smooth $\P X\subset \P V$ the projective dual of
$\P X$ is simply the variety of projective hyperplanes tangent to
$\P X$.

Our objective in this section is to find a formula for the
$T$-equivariant rational cohomology class of the dual $\check X$ of
$X\subset V$, under certain conditions ($\check X$ is the cone of
$\P\check X$). The result will be a cohomology class in $H^*(BT)$
which we identify with the symmetric algebra of the character group
of $T$. Since we use localization, our formulas will formally live
in the fraction field of $H^*(BT)$, but part of the statements will
be that in the outcome the denominators cancel.

Let $\check{V}$ denote the dual vector space of $V$ and consider the
conormal vector bundle
$$\check{N}=\{(x,\Lambda)\in \P X\times \check{V}\ |\ \P\Lambda
\text{ is tangent to } \P X \text{ at } x\},$$ and the diagram
$$\xymatrix{ \check{N} \ar@/^2pc/[rrr]^{\phi} \ar@{^{(}->}@<-3pt>[rr]^{i}  \ar[rd] &  & \P X\times
\check{V} \ar[ld]
\ar[r]^{\pi} & \check{V} \\
 & \P X &  & \check{X} \ar@{^{(}->}[u] } $$
(cf. Proposition~\ref{31}).

One can check that the assumption that $\P \check{X}$ is a
hypersurface yields that the map $\phi$ is a resolution of
$\check{X}$ (see \cite[thm 1.10]{tevelev}). Hence from (\ref{we}) we
obtain the following statement.


\begin{proposition} Assume that the fixed point set $F(\P X)$ is finite.
Then the equivariant cohomology class represented by the cone
$\check X$  is
\[[\check X]=\sum_{f\in F(\P X)}\frac{e(\check V)}{e(T_f\P X)e(\check N_f\P X)}.\]
\end{proposition}

Let $\W(A)$ denote the set of weights of the $T$-module $A$ counted
with multiplicity.  For a fixed point $f\in F(\P X)$ we have
\[\W(\check V)=\{-\omega (f)\} \cup \{-\omega (f)-\beta| \beta\in \W (T_f \P X)\} \cup \W(\check{N}_f\P X),\]
where $\omega (f)$ denotes the weight corresponding to the fixed
point $f$. Notice the asymmetrical role of $\check N_f\P X$ and
$T_f\P X$. The vector space $\check N_f\P X$ is a subspace of
$\check V$ but $T_f\P X\subset T_f\P V\iso \hom(L_f, V\ominus L_f)$
where $L_f\leq V$ is the eigenline corresponding to the fixed point
$f$. This isomorphism explains the shift
by $-\omega (f)$ in the formula. Therefore we obtain 
\begin{theorem}
 \label{main} Suppose that  the torus $T$ acts on the complex vector space $V$ linearly and $X\subset V$ is a $T$-invariant cone. Assume moreover that $\P X\subset\P  V$ is smooth with finitely many $T$-fixed points and the projective dual $\P \check X$ is a hypersurface. Then the equivariant cohomology class represented by the cone $\check X$ of the dual $\P\check X$ is
     $$[\check X]=-\sum_{f\in F(\P X)}\omega (f)\prod_{\beta\in \W(T_f\P X)}\frac{\beta+\omega (f)}{-\beta}.$$
\end{theorem}

Notice that the fixed points $f$ of the $T$-action on $\P X$
correspond to eigenlines $L_f\leq V$ and the weight $\omega (f)$ of
this line is canonically identified with the cohomology class
$e(L_f)$.

\begin{remark} \label{posdef}
   The difference of the dimensions of a hypersurface and the
   variety $\P \check X$ is called the {\em defect} of $\P X$.
   Hence Theorem \ref{main} deals with the defect 0 case.
   It is customary to define the {\em cohomology class} of $\check X$ and the {\em degree}
   of $\P  \check X$ to be 0 if $\P X$ has positive defect. Using this convention Theorem~\ref{main}
   remains valid without the condition on the defect. Indeed, by Remark~\ref{intzero} the right hand side, which is equal to the pushforward of 1, is automatically zero if the image of $\pi_2$ has smaller dimension.
\end{remark}

\begin{remark}
Similar argument yields the cohomology class of $\P X$ itself: In
Proposition~\ref{31} we choose $S$ to be the tautological bundle
over $\P X$. As a result we obtain a nontrivial special case of the
{\bf Duistermaat-Heckman formula}: {\sl If the set $F(\P X)$ of
fixed points on $\P X$ is finite, then
$$ [X]=e(V)\sum_{f\in F(\P X)}(\omega (f)\cdot e(T_f\P X))^{-1}.$$}

\end{remark}

\section{Cohomology and degree formulas for the discriminant}\label{degreeformulas}

In this section we apply Theorem \ref{main} to obtain formulas for
the equivariant class and the degree of the discriminants of
irreducible representations.

Let \( \rho:G\to GL(V) \) be an irreducible representation of the
complex connected reductive Lie group $G$ on the complex vector
space $V$, and let $\P X\subset \P V$ be the (closed, smooth) orbit
of $[v]$ where $v\in V$ is a vector corresponding to the highest
weight $\lambda$.
Let $D_\lambda \subset \check{V}$ and $\p D_{\lambda}\subset \P
\check{V}$ be the duals of $X$ and $\P X$, they are called the
``discriminants'' of the representation of highest weight $\lambda$.
The discriminant is ``usually'' a hypersurface, a complete list of
representations of semisimple Lie groups for which the discriminant
is not a hypersurfare (i.e. the defect is positive, c.f.
Remark~\ref{posdef}) can be found in \cite{km}, see also
\cite[Th.9.21]{tevelev}.

\begin{remark}
  For an irreducible representation of a reductive group $G$ on $V$ the action of the center of $G$ is trivial on $\P V$, that is,
  only the semisimple part of $G$ acts on $\P V$. Hence, seemingly it is enough to state the degree formula only for semisimple groups. However, we will calculate the degree of the discriminant from its
  equivariant Poicar\'e dual cohomology class, which is an element in $H^2_G($point$)$. For semisimple groups this
  cohomology group is 0. Therefore, if $G$ is semisimple, we replace $G$ with $G\times GL(1)$---and let $GL(1)$ act on $V$ by multiplication.
  Consequently we are forced to state our theorems for reductive
groups.
\end{remark}

Let $W$ be the Weyl group, let $R(G)$ be the set of roots of $G$,
and $R^-(G)$ the set of negative roots. Let $\g_\beta$ be the root
space corresponding to the root $\beta$. Let $H_\beta$ be the unique
element in $[\g_\beta,\g_{-\beta}]$ with $\beta(H_\beta)=2$.


\begin{theorem} \label{short} {\bf Main Formula---Short Version.}
With $G$, $W$, $V$, $\lambda$, $D_\lambda$, $R^-(G)$ as above, the
equivariant Poincar\'e dual of $D_\lambda$ in $\check{V}$ is
$$[D_\lambda]=-\sum_{\mu \in W\lambda} \mu \prod_{\beta\in T_\mu}
\frac{\mu+\beta}{-\beta},$$ where
$$T_\lambda=\{\beta \in R^-(G)\ |\  \langle  H_\beta,\lambda \rangle<0\},$$
and $T_{w\lambda}=wT_\lambda$ for $w\in W$. Here we used the
convention of Remark~\ref{posdef}, i.e. the class of $[D_\lambda]$
is defined to be 0 if $D_\lambda$ is not a hypersurface.
\end{theorem}

Let us remark that $\{\beta \in R^-(G)| B(\beta,\lambda)<0\}$ is an
equivalent description of $T_\lambda$ (where $B$ is the Killing
form). The proof of Theorem \ref{short} is based on the following
two standard lemmas:
\begin{lemma}\label{fixedpoints}
The fixed point set $F(\P X)$ of the maximal torus $T\subset G$ is
equal to the orbit of  $[v]\in\P V$ for the action of the Weyl group
$W$.
\end{lemma}
Notice that the  Weyl group $W=N_G(T)/T$ indeed acts on $\P X$ since
$T$ fixes $[v]\in \P V$.
\begin{proof}
It is enough to show that if $[v]\in \P V$ and $g[v]$ are both fixed
points of $T$ and $v$ is a maximal weight vector, then there exists
a $\beta\in N_G(T)$ such that $g[v]=\beta[v]$.

Let $G_{[v]}$ be the stabilizer of $[v]$. Then, by the assumption,
$T$ and $g^{-1}Tg$ are contained in  $G_{[v]}$. These are maximal
tori in  $G_{[v]}$, so there is a $p\in  G_{[v]}$ such that
$g^{-1}Tg=p^{-1}Tp$ (see e.g. \cite[p. 263]{borel-lag}). Then
$\beta=gp^{-1}\in N_G(T)$ and $g[v]=\beta[v]$.
\end{proof}
\begin{lemma} \label{tangent} The weights of the tangent space $T_f(\P X)$ as a $T$-space are
$$T_f=\{\beta \in R^-(G)\ |\  \langle H_\beta,\omega (f) \rangle <0\}$$
for any  $f\in F(\P X)$.
\end{lemma}

For semisimple $G$ the proof can be found in \cite{fulton-harris}
and a more detailed version in \cite[p.36]{berczi}.
The formula extends to the reductive case without change.

\begin{proof}[Proof of Theorem~\ref{short}] It is enough to apply Theorem~\ref{main} to our situation.
Lemma~\ref{fixedpoints} determines the (finitely many) fixed points
and  Lemma~\ref{tangent} gives that
$$T_\lambda=\{\beta \in R^-(G)\ |\  \langle  H_\beta,\lambda \rangle<0\}.$$
The weights at other fixed points are obtained by applying the
appropriate element of the Weyl group.
\end{proof}

For the Lie group $G=GL(n)$ the maximal torus can be identified with
the subgroup of diagonal matrices $diag(z_1,\ldots,z_n)$, $|z_i|=1$,
which is the product of $n$ copies of $S^1$'s. Let $L_i$ be the
character of this torus, which is the identity on the $i$'th $S^1$
factor, and constant 1 on the others. The Weyl group, the symmetric
group on $n$ letters, permutes the $L_i$'s. The Weyl-invariant
subring of the symmetric algebra $\Q[L_1,\ldots,L_n]$ is identified
with $H^*(BGL(n))$. Under this identification, the $i$'th elementary
symmetric polynomial of the $L_i$'s is the $i$'th Chern class
(`splitting lemma'). Hence---as usual---we call the $L_i$'s the {\em
Chern roots} of $GL(n)$.

\begin{example} Consider the dual of $Gr_3(\C^8)$ in its Pl\"ucker
embedding. This is the discriminant of the representation of $GL(8)$
with highest weight $L_1+L_2+L_3$. Let ${n \choose k}$ stand for the
set of $k$-element subsets of $\{1,\ldots,n\}$.
  By Theorem
\ref{short} the class of the discriminant is
$$-\sum_{S\in {8\choose 3}} (L_{s_1}+L_{s_2}+L_{s_3})\cdot
\prod_{i\in S} \prod_{j\not\in S}
\frac{L_{s_1}+L_{s_2}+L_{s_3}+L_j-L_i}{(L_i-L_j)}.$$ This is a
58-term sum. However, this class is in $H^2(BGL(8))$, hence we know
that it must have the form of $v\cdot c_1=v\cdot\sum_{i=1}^8 L_i$.
Therefore we only need to determine the number $v$. A well chosen
substitution will kill most of the terms, making the calculation of
$v$ easier. E.g.~substitute $L_i=i-13/3$, then all terms---labeled
by 3-element subsets $S$ of $\{1,\ldots,8\}$---are zero, except for
the last one corresponding to $S=\{6,7,8\}$. For the last term we
obtain $-8$. Hence, $-8=v\cdot (1-13/3+2-13/3+\ldots+8-13/3)$, which
gives the value of $v=-6$. Thus we obtain the equivariant class of
the dual of $Gr_3(\C^8)$ to be $-6(L_1+\ldots+L_8)$, and in turn its
degree as $-6(-1/3-\ldots-1/3)=16$ (by Proposition~\ref{deg}).
\end{example}

Although similarly lucky substitutions cannot be expected in
general, the $L_i=i$ substitution yields a formula for the degree of
the dual variety of the Grassmannian $Gr_k(\C^n)$ in its Pl\"ucker
embedding:
\begin{equation}\label{grass}
\deg(\check{Gr}_k(\C^n))= \frac{2k}{n+1}\sum_{S\in {n \choose k}}
l(S) \prod_{i\in S, j\not\in S} \frac{l(S)+j-i}{i-j},\end{equation}
 where $l(S)=\sum_{s\in S} s$.

\medskip

To study how the degree depends on the highest weight $\lambda$ for
a fixed group $G$ we introduce another expression for the degree
where the sum is over all elements of the Weyl group instead of the
orbit $W\lambda$.

For a dominant weight $\lambda$, let $O_\lambda=R^-(G)\setminus
T_\lambda=\{\beta\in R^-(G)\ |\  \langle
H_\beta,\lambda\rangle=0\}$, and let the sign $\varepsilon(\lambda)$
of $\lambda$ be $(-1)^{|O_\lambda|}$. Let $W_\lambda\leq W$ be the
stabilizer subgroup of $\lambda$.

\begin{theorem} \label{long} {\bf Main Formula---Symmetric Version.}
Under the conditions of Theorem~\ref{short} we have
$$[D_\lambda]=-\frac{\varepsilon(\lambda)}{|W_{\lambda}|}
\sum_{w\in W} w\Big( \lambda \prod_{\beta\in R^-(G)}
\frac{\lambda+\beta}{-\beta}\Big).$$
\end{theorem}

\begin{proof} Let $w_1,\ldots,w_m$ be left coset representatives
of $W_\lambda\leq W$, i.e. the disjoint union of the
$w_iW_\lambda$'s is $W$. Then
$$
-\frac{\varepsilon(\lambda)}{|W_{\lambda}|} \sum_{w\in W} w\left(
\lambda \prod_{\beta\in R^-(G)}
\frac{\lambda+\beta}{-\beta}\right)=$$
\begin{align}
&= -\frac{\varepsilon(\lambda)}{|W_{\lambda}|}\sum_{i=1}^m
w_i\left(\sum_{w\in W_\lambda} w\left( \lambda \prod_{\beta \in
T_\lambda} \frac{\lambda+\beta}{-\beta}\prod_{\beta \in O_\lambda}
\frac{\lambda+\beta}{-\beta} \right)\right) \\
&= \label{tovabb}
-\frac{\varepsilon(\lambda)}{|W_{\lambda}|}\sum_{i=1}^m w_i\Big(
\lambda \prod_{\beta \in T_\lambda}
\frac{\lambda+\beta}{-\beta}\sum_{w\in W_\lambda} w\big(\prod_{\beta
\in O_\lambda} \frac{\lambda+\beta}{-\beta} \big)\Big),
\end{align}
since $wT_\lambda=T_\lambda$ (but $wO_\lambda$ is not necessarily
equal to $O_\lambda$). Now we need the following lemma.
\begin{lemma}
$$\sum_{w\in W_\lambda} w\big( \prod_{\beta\in O_\lambda}
\frac{\lambda+\beta}{-\beta}\big)=\varepsilon(\lambda)
|W_\lambda|.$$
\end{lemma}
\begin{proof} The polynomial
$$P(x)= \sum_{w\in W_\lambda} \prod_{\beta \in O_\lambda}
\frac{x+w(\beta)}{-w(\beta)} \in \Z(L_i)[x]$$ is
$W_\lambda$-invariant, and has the form
$$\frac{q_k+q_{k-1}x+\ldots+q_0x^k}{\prod_{\beta \in
O_\lambda} \beta}$$ where $k=|O_\lambda|$ and $q_i$ is a degree $i$
polynomial on the orthogonal complement of $\lambda$. This means
that the numerator has to be anti-symmetric under $W_\lambda$ (which
is itself a Weyl group of a root system), hence it must have degree
at least the number of positive roots. That is, all $q_i$, $i<k$
must vanish. This means that $P(x)$ is independent of $x$, i.e.
$P(\lambda)=P(0)=\varepsilon(\lambda)|W_\lambda|$, as required.
\end{proof}
Using this Lemma, formula (\ref{tovabb}) is further equal to
$$-\sum_{i=1}^m w_i \Big( \lambda \prod_{\beta\in T_\lambda}
\frac{\lambda+\beta}{-\beta}\Big),$$ which completes the proof of
Theorem \ref{long}.
\end{proof}

The advantage of the Short Version (Theorem \ref{short}) is that for
certain $\lambda$'s ($\lambda$'s on the boundary of the
Weyl-chamber) the occurring products have only few factors, while
the advantage of the Symmetric Version (Theorem \ref{long}) is that
it gives a unified formula for all $\lambda$'s of a fixed group. Now
we will expand this latter observation.

Let $G$ be semisimple and consider the representation with highest
weight $\lambda$. Extend this action to an action of $G\times GL(1)$
with $GL(1)$ acting by scalar multiplication. Denoting the first
Chern class of $GL(1)$ by $u$ (that is $H^*(BU(1))=\Q[u]$) we obtain
that
\begin{equation} \label{tpu}
[D_\lambda]=-\frac{\varepsilon(\lambda)}{|W_{\lambda}|} \sum_{w\in
W} w\Big( (\lambda+u) \prod_{\beta\in R^-(G)}
\frac{\lambda+u+\beta}{-\beta}\Big). \end{equation}
 Proposition \ref{deg} then turns this formula to a degree formula for the
discriminant, by substituting $u=-1$:
\begin{equation}\deg(D_\lambda)=-\frac{\varepsilon(\lambda)}{|W_{\lambda}|}
\sum_{w\in W} w\Big( (\lambda-1) \prod_{\beta\in R^-(G)}
\frac{\lambda+\beta-1}{-\beta}\Big).
\end{equation}
Part of the statement is that this formula is a constant, i.e.
expression (\ref{tpu}) is equal to a constant times $u$ (although
this can also be deduced from the fact that $H^2(BG)=0$ for
semisimple groups).

\begin{corollary} \label{fehpol} Let $G$ be a semisimple Lie group,
and $\h$ the corresponding Cartan subalgebra. There exists a
polynomial $F_G:\h^*\to \Z$, with degree equal the number of
positive roots of $G$, such that
$$\deg(D_\lambda)=\frac{\varepsilon(\lambda)}{|W_{\lambda}|}
F_G(\lambda)$$ if $\lambda$ is a dominant weight. The polynomial
$F_G(\lambda)$ vanishes if and only if $D_\lambda$ is not a
hypersurface.
\end{corollary}

\begin{proof} We have
\begin{equation} \label{feherpoly}
F_G=-\sum_{w\in W} w\Big( (\lambda-1) \prod_{\beta\in R^-(G)}
\frac{\lambda+\beta-1}{-\beta}\Big), \end{equation} and the last
statement follows from Remark \ref{posdef}.
\end{proof}

\begin{remark} The polynomial dependence of $\deg(D_\lambda)$ for
{\em regular} weights $\lambda$ (hence $\varepsilon(\lambda)=1$,
$W_\lambda=\{1\}$), as well as a formula for a special value of the
polynomial (the value at the sum of the fundamental weights) is
given in \cite{weyman}. Since $\deg(D_\lambda)$ is always
non-negative the value of $\varepsilon(\lambda)$ is determined by
the sign of $F_G(\lambda)$. The positive defect cases are known but
this corollary provides an alternative and uniform way to find them.
\end{remark}

\begin{example}\label{akbek}
A choice of simple roots $\alpha_1,\alpha_2,\ldots,\alpha_r$ of $G$
determines the fundamental weights $\omega_1,\omega_2,\ldots,
\omega_r$ by
$B(\omega_i,\frac{2\alpha_j}{B(\alpha_j,\alpha_j)})=\delta_{i,j}$
where $B(.,.)$ denotes the Killing form. We follow the convention of
De~Concini and Weyman \cite{weyman} by writing $F_G$ in the basis of
fundamental weights (i.e. $y_1\omega_1+\ldots+y_r\omega_r\mapsto
F_G(y_1,\ldots,y_r)$) and substituting $y_i=x_i+1$. The advantage of
this substitution is that in this way, according to  \cite{weyman},
all the coefficients of the polynomial $F_G$ become non-negative.
Formula (\ref{feherpoly}) gives the following polynomials for all
semisimple Lie groups of rank at most 2 and for some of rank 3. [In
these examples our convention for simple roots agrees with the one
in the {\it coxeter/weyl} Maple package {\sl
www.math.lsa.umich.edu/\~{}jrs/maple.html} written by J.
Stembridge.]

\vspace{4mm}

\noindent  $ {\mathbf{ A_1}} $,\ \ \ \ $\alpha_1=L_2-L_1$:

\vspace{1mm}

$F=2x_1$.

\vspace{2mm}

\noindent  ${\mathbf{ A_1+ A_1}}$,\ \ \ \
$\alpha=(L_2-L_1,L_2'-L_1')$:

\vspace{1mm}

$F=6x_1x_2+2x_2+2x_1+2$.

\vspace{2mm}

\noindent  ${\mathbf{ A_2}}$,\ \  \ \ $\alpha=(L_2-L_1, L_3-L_2)$:

\vspace{1mm}

$F=6(x_1+x_2+1)(2x_1x_2+x_1+x_2+1)$ \cite[ex 7.18]{tevelev}.

\vspace{2mm}

\noindent  ${\mathbf{ B_2}}$,\ \ \ \  $\alpha=(L_1, L_2-L_1)$,

\vspace{1mm}

$F=20(2x_2^3x_1+3x_2^2x_1^2+x_2x_1^3)+12(2x_2^3+12x_2^2x_1+
11x_2x_1^2+x_1^3)+24(3x_2^2+7x_2x_1+2x_1^2)+8(9x_2+8x_1)+24$.

\vspace{2mm}

\noindent ${\mathbf{ G_2}}$,\ \ \ \  $\alpha=(L_2-L_1,
L_1-2L_2+L_3)$:

\vspace{1mm}

$F=42(18x_2^5x_1+45x_2^4x_1^2+40x_2^3x_1^3+15x_2^2x_1^4+2x_2x_1^5)
+60(9x_2^5+90x_2^4x_1+150x_2^3x_1^2+90x_2^2x_1^3+20x_2x_1^4+x_1^5)
+110(27x_2^4+132x_2^3x_1+144x_2^2x_1^2+52x_2x_1^3+5x_1^4)+
8(822x_2^3+2349x_2^2x_1+1527x_2x_1^2+248x_1^3)+
6(60x_2^2+1972x_2x_1+579x_1^2)+4(1025x_2+727x_1)+916$.

\vspace{2mm}

\noindent  ${\mathbf{ A_1+A_1+A_1}}$,\ \ \ \
$\alpha=(L_2-L_1,L_2'-L_1',L_2''-L_1'')$:

\vspace{1mm}

$F=24x_2x_1x_3+12(x_2x_1+x_2x_3+x_1x_3)+8(x_2+x_1+x_3)+4$.

\vspace{2mm}

\noindent ${\mathbf{ A_1+A_2}}$,\ \ \ \  $\alpha=(L_2-L_1,L_2'-L_1',
L_3'-L_2')$:

\vspace{1mm}

$F=60x_2^2x_1x_3+60x_2x_1x_3^2+36x_2^2x_1+36x_2^2x_3+144x_2x_1x_3+
36x_2x_3^2+36x_1x_3^2+24x_2^2+72x_2x_1+96x_2x_3+72x_1x_3+24x_3^2+48x_2+36x_1+48x_3+24$.

\vspace{2mm}

\noindent ${\mathbf{ A_3}}$,\ \ \ \  $\alpha=(L_2-L_1, L_3-L_2,
L_4-L_3)$:

\vspace{1mm}

$F=420x_1x_3x_2(x_1+x_2)(x_2+x_3)(x_2+x_1+x_3)+
300(x_1^3x_2^2+4x_1^3x_2x_3+x_1^3x_3^2+2x_1^2x_2^3+15x_1^2x_2^2x_3+
12x_1^2x_2x_3^2+x_1^2x_3^3+x_1x_2^4+12x_1x_2^3x_3+15x_1x_2^2x_3^2+
4x_1x_2x_3^3+x_2^4x_3+2x_2^3x_3^2+x_2^2x_3^3)+
220(3x_1^3x_2+3x_1^3x_3+12x_1^2x_2^2+33x_1^2x_2x_3+9x_1^2x_3^2+
10x_1x_2^3+48x_1x_2^2x_3+33x_1x_2x_3^2+3x_1x_3^3+x_2^4+10x_2^3x_3+12x_2^2x_3^2+3x_2x_3^3)+
8(42x_1^3+453x_1^2x_2+411x_1^2x_3+699x_1x_2^2+1566x_1x_2x_3+
411x_1x_3^2+164x_2^3+699x_2^2x_3+453x_2x_3^2+42x_3^3)+
12(126x_1^2+491x_1x_2+407x_1x_3+239x_2^2+491x_2x_3+126x_3^2)+
16(133x_1+169x_2+133x_3)+904 $.

\vspace{2mm}

 A Maple computer program computing $F$ for any semisimple Lie
group is available at
 www.unc.edu/\-\~{}rimanyi\-/progs\-/feherpolinom.mw.

\vspace{2mm}

 Straightforward calculation gives that $F_G$ for
$G=A_1+\ldots+ A_1$ ($n$ times) is
\[F_{nA_1}(y_1,\ldots,y_n)=\sum_{k=0}^{n} (-2)^{n-k}(k+1)!\, \sigma_k(y_1,\ldots,y_n),\]
where $\sigma_i$ is the $i$'th elementary symmetric polynomial.
This can also be derived from \cite[Th.2.5,Ch.13]{gkz}.

\end{example}

\section{Combinatorics of the degree formulas for $GL(n)$} \label{combinat}
Our main formulas, Theorems \ref{short} and \ref{long}, can be
encoded using standard notions from the combinatorics of symmetric
functions. In Sections \ref{combinat}-\ref{andras} we assume that
$G=GL(n)$, and that the simple roots are $L_i-L_{i+1}$. Then
$R^-(G)=\{L_i-L_j\,:\, i>j\}$. Irreducible representations
correspond to the weights $\lambda=\sum_{i=1}^na_iL_i$ with $a_1\geq
a_2\geq \cdots \geq a_n$. For such a weight one has
$T_\lambda=\{L_i-L_j\,:\, i>j,\ a_i<a_j\}$.
\subsection{Symmetrizer operators} \label{symmet}
Let
$$\lambda^+=\lambda\prod_{\beta\in R^-(G)} (\lambda+\beta),\qquad\hbox{and}\qquad
\Delta=\prod_{1\leq i < j\leq n}(L_i-L_j),$$ and recall the
definition of the Jacobi-symmetrizer (\cite{lascouxbook}) of a
polynomial $f(L_1,\ldots,L_n)$:
$$J(f)(L_1,\ldots,L_n)=\frac{1}{\Delta}\sum_{w\in S_n} \sgn(w) f(L_{w(1)},\ldots,L_{w(n)}),$$
 where $\sgn(w)$ is the sign of the permutation $w$, i.e. $(-1)$ raised to the
 power of the number of transpositions in $w$. Then for
$\lambda=\lambda_1L_1+\ldots+\lambda_nL_n$ with $\lambda_1\geq
\lambda_2\geq\ldots\geq \lambda_n$, from Theorem \ref{long} we
obtain
\begin{equation} \label{jacobi}
[D_\lambda]=-\frac{\varepsilon(\lambda)}{|W_\lambda|}J\Big(\lambda^+\Big),
\end{equation}
and $\deg D_\lambda$ is obtained by multiplying this  by the factor
\begin{equation}\label{degdlambda}
 \frac{-n}{|\lambda|\cdot \sigma_1}, \ \ \ \mbox{where} \ \
|\lambda|:= \sum \lambda_i  \ \ \mbox{and} \ \ \sigma_1:=\sum
L_i.\end{equation}

\def\d{\partial}

The Jacobi symmetrizer is a special case of the divided difference
operators $\d_w$ \cite[Ch.7]{lascouxbook}, corresponding to the
maximal permutation $w_0=[n,n-1,\ldots,2,1]$ (i.e. $w_0(i)=n+1-i$).
Shorter divided difference operators also turn up for certain
representations. We will illustrate this with the case of the
Pl\"ucker embedding of Grassmann varieties, i.e. $G=GL(n)$,
$\lambda=L_1+\ldots+L_k$. Let $[k|n-k]\in S_n$ be the permutation
$[n-k+1,n-k+2,\ldots,n,1,2,\ldots,n-k]$. Theorem \ref{short} gives
$$[D_{L_1+\ldots+L_k}]=-\d_{[k|n-k]}\Big( (L_1+\ldots+L_k)
\prod_{i=1}^{k} \prod_{j=k+1}^{n} (L_1+\ldots+L_k+L_j-L_i)\Big),$$
and the degree of the discriminant of $Gr_k(\C^n)$ in its Pl\"ucker
embedding is obtained by multiplying this with $-n/(k \sigma_1)$.

\subsection{Scalar product} Now we show how to use the scalar product
on function spaces to encode the formula of Theorem \ref{long} in
the case of $G=GL(n)$. Following \cite{macdonald}, we define the
scalar product of the polynomials $f,g$ in $n$ variables
$L_1,\ldots,L_n$ as
$$\langle f,g\rangle=\frac{1}{n!}
\Big[f\bar{g}\prod_{i\not=j}\big(1-\frac{L_i}{L_j}\big)\Big]_1,$$
where $\bar{g}(L_1,\ldots,L_n)=g(1/L_1,\ldots,1/L_n)$ and $[h]_1$ is
the constant term (i.e. the coefficient of 1) of the Laurent
polynomial $h\in \Z[L_1^{\pm},\ldots,L_n^{\pm}]$. The Jacobi
symmetrizer is basically a projection, thus for a degree ${n \choose
2}+1$ polynomial $f$ we have
$$\frac{1}{n!}\sum_w \sgn(w)
f(L_{w(1)},\ldots,L_{w(n)})= \frac{\langle
f,\sigma_1\Delta\rangle}{\langle \sigma_1\Delta,\sigma_1
\Delta\rangle}\sigma_1\Delta.$$

Here $\langle \sigma_1\Delta,\sigma_1\Delta\rangle$ can be
calculated to be $n(2n-3)!!$, where $(2k+1)!!$ denotes the
semifactorial $1\cdot3\cdots(2k-1)\cdot(2k+1)$, hence we obtain
$$[D_\lambda]=\frac{-\varepsilon(\lambda)n!}{|W_\lambda|\cdot n(2n-3)!!} \langle\lambda^+, \sigma_1\Delta\rangle \sigma_1,$$
and hence the following form of our Main Formula:

\begin{theorem} The degree of the discriminant of the irreducible
representation of $GL(n)$ with highest weight $\lambda$ is
$$\deg D_\lambda=\frac{\varepsilon(\lambda)}{|\lambda|\cdot |W_\lambda|}
\frac{n!}{(2n-3)!!}
 \langle\lambda^+, \sigma_1\Delta\rangle.$$
\end{theorem}

\subsection{Permanent} \label{perma}
For $(\nu_1,\ldots,\nu_n)\in \N^n$ and $w\in S_n$ let
$w(\nu)=(\nu_{w(1)},$ $\ldots,$ $\nu_{w(n)})$ and $L^\nu$ will
denote the monomial $L_1^{\nu_1}\ldots L_n^{\nu_n}$.

\begin{lemma} \label{j}
 Let $\mu=(n,n-2,n-3,\ldots,2,1,0)\in \N^n$. If $\sum \nu_i={n\choose 2}+1$ then we have
$$J(L^\nu)=\begin{cases}\sgn(w)\sigma_1 & \hbox{if}\ \nu_i=\mu_{w(i)} \\
0 & \hbox{otherwise} \end{cases}$$
\end{lemma}

\begin{proof} If $\nu_i=\nu_j$ then the terms of $J(L^\nu)$ turn
up in cancelling pairs, hence $J(L^\nu)=0$. This leaves only $\nu=$
permutations of $\mu$ for possible non-zero $J$-value. Direct check
shows $J(L^\mu)= \sigma_1$ (c.f. the well known identity
$J\left(L^{(n-1,n-2,\ldots,2,1,0)}\right)=1$).
\end{proof}

The coefficient of a monomial $L^\nu$ in a polynomial $f$ will be
denoted by $c(f,L^\nu)$. Formula (\ref{jacobi}) and Lemma \ref{j}
yield that the class and degree of $D_\lambda$ can be computed by
counting coefficients.

\begin{theorem} \label{coef}
The equivariant class of $D_\lambda$ is
$$[D_\lambda]=-\frac{\varepsilon(\lambda)}{|W_\lambda|} \sum_{w\in
S_n} \sgn(w) c(\lambda^+,L^{w(\mu)}) \cdot \sigma_1,$$ and the
degree of $D_\lambda$ is obtained by multiplication by
$-n/(|\lambda|\sigma_1)$.
\end{theorem}

Similar sums will appear later, hence we define the $\nu$-permanent
of a polynomial $f$ as
$$\sum_{w\in S_n} \sgn(w) c(f,L^{w(\nu)}),$$
and denote it by $P(f,\nu)$.

The name {\em permanent} is justified by the following observation.
Let $\nu=(\nu_1,\ldots,\nu_n)$ be a partition. If the polynomial $f$
is the product of $|\nu|$ linear factors $\sum_{j=1}^n
a^{(i)}_jL_j$, $i=1,\ldots,|\nu|$, then $P(f,\nu)$ can be computed
from the $n\times |\nu|$ matrix $(a^{(i)}_j)$ as follows. For a
permutation $w\in S_n$ we choose $\nu_1$ entries from row $\nu(1)$
(i.e. $j=\nu(1)$), $\nu_2$ entries from row $\nu(2)$, etc, such a
way that no chosen entries are in the same column (hence they form a
complete `rook arrangement'). The product of the chosen entries with
sign $\sgn(w)$ will be a term in $P(f,\nu)$, and we take the sum for
all $w\in S_{n}$ and all choices. For example the $(2,1)$-permanent
of the product $(aL_1+a'L_2)(bL_1+b'L_2)(cL_1+c'L_2)$, considering
the matrix $\begin{pmatrix} a & b & c \\ a' & b' & c'
\end{pmatrix}$, is
\[ abc'+ab'c+a'bc-a'b'c-a'bc'-ab'c'.\]

\begin{theorem}\label{boole54} (Boole's formula, \cite[7.1]{tevelev})
If $\lambda=aL_1$, then the degree of $D_\lambda$ is $n(a-1)^{n-1}$.
\end{theorem}

\begin{proof} We have
$$\lambda^+=aL_1\prod_{i=2}^n \left((a-1)L_1+L_i\right)\cdot\mspace{-12mu} \prod_{2\leq i<j \leq n}(L_j-L_i+aL_1),$$
and let
$$\lambda^{++}=aL_1\prod_{i=2}^n ((a-1)L_1+L_i)\cdot \mspace{-12mu}\prod_{2\leq i<j \leq n}(L_j-L_i).\phantom{+aL_1+}$$

Observe that $P(\lambda^+,\mu)=P(\lambda^{++},\mu)$, since in the
difference each term comes once with a positive, once with a
negative sign. In $P(\lambda^{++},\mu)$ only $L_1$ has degree $\geq
n$ (namely $n$), hence we have
$$P(\lambda^{++},\mu)=a(a-1)^{n-1}P\left(\prod_{2\leq i<j\leq n}(L_i-L_j),(n-2,n-3,\ldots,2,1)\right)$$
which is further equal to  $a(a-1)^{n-1}(-1)^{{n-1 \choose 2}}
(n-1)!$ and Theorem \ref{coef} gives the result.
\end{proof}

\subsection{Hyperdeterminants}  The discriminant of the standard action of the
product group $\prod_{u=1}^k GL(n_u)$ on
$\bigotimes_{u=1}^k\C^{n_u}$ is called {\em hyperdeterminant}
because it generalizes the case of the determinant for $k=2$. In
other words, the hyperdeterminant is the dual of the Segre embedding
$\prod_{u=1}^k\P(\C^{n_u})\to \P(\bigotimes_{u=1}^k\C^{n_u})$.
Gelfand, Kapranov and Zelevinsky \cite{gkz} give a concrete
description of the degree of the hyperdeterminants by giving a
generator function. In the so-called boundary case when
$n_k-1=\sum_{u=1}^{k-1} (n_u-1)$ they get a closed formula. We show
how to prove this using $\nu$-permanents.

When considering representations of $\prod_{u=1}^k GL(n_u)$
($n_1\leq n_2\leq \ldots \leq n_k$) we will need $k$ sets of
variables (the $k$ sets of Chern roots), we will call them
$L_{(u),i}$, $u=1,\ldots,k$, $i=1,\ldots,n_u$. The hyperdeterminant
is the discriminant corresponding to the representation with highest
weight $\lambda=\sum_{u=1}^k L_{(u),1}$. If $\nu_{(u)}\in \N^{k_u}$,
then a polynomial $f$ in these variables has a
$\nu=(\nu_{(1)},\ldots,\nu_{(k)})$-permanent, defined as
$$P(f,\nu)=\sum_{w^{(1)}\in S_{n_1}}\ldots\sum_{w^{(k)}\in S_{n_k}}
\big(\prod_{u=1}^k \sgn(w^{(u)})\big) c(f,\prod_u
L_{(u)}^{w^{(u)}(\nu_{(u)})}).$$

\begin{theorem}[{\cite[14.2.B]{gkz}}]  If $n_k-1>\sum_{u=1}^{k-1} (n_u-1)$ then the
discriminant is not a hypersurface. If $n_k-1=\sum_{u=1}^{k-1}
(n_u-1)$ then its degree is $n_k!/\prod_{u=1}^{k-1}(n_u-1)!$.
\end{theorem}

\begin{proof} Let $\mu^{(u)}=(\mu_{(1)},\ldots,\mu_{(k)})$, where $\mu_{(v)}=(n_v-1,n_v-2,\ldots,1,0)$ for $v\not=u$, and
$\mu_{(u)}=(n_u,n_u-2, n_u-3,\ldots,2,1,0)$. A straightforward
generalization of the $k=1$ case gives
\begin{equation} \label{hyp}
[D_\lambda]=-\frac{(-1)^{\sum {n_u-1 \choose
2}}}{\prod_{u=1}^k(n_u-1)!}\sum_{u=1}^k
\Big(P(\lambda^+,\mu^{(u)})\sum_{i=1}^{n_u}L_{(u),i}\Big),
\end{equation}
where
$$\lambda^+=\lambda \left(\prod_{u=1}^k
\prod_{i=2}^{n_u}(L_{(u),i}-L_{(u),1}+\lambda)\right)\left(\prod_{u=1}^k
\mspace{8mu} \prod_{2\leq i<j\leq n_u}
(L_{(u),j}-L_{(u),i}+\lambda)\right).$$ Just like in the proof of
the Boole formula above, we can change $\lambda^+$ to
$$\lambda^{++}=\lambda \left(\prod_{u=1}^k
\prod_{i=2}^{n_u}(L_{(u),i}-L_{(u),1}+\lambda)\right)\left(\prod_{u=1}^k\mspace{8mu}
\prod_{2\leq i<j\leq n_u} (L_{(u),j}-L_{(u),i})\right)\phantom{
++\lambda}$$ without changing the value of the permanent. Observe
that the highest power of an $L_{(k),i}$ variable in $\lambda^{++}$
is $\sum_{u=1}^{k-1} (n_u-1)$, thus the last term in (\ref{hyp}) is
0 if $n_k-1>\sum_{u=1}^{k-1} (n_u-1)$. This proves the first
statement of the theorem.

If we have $n_k-1=\sum_{u=1}^{k-1} (n_u-1)$, then considering the
power of $L_{(k),1}$ again, we get that

\begin{multline*}
  P(\lambda^{++},\mu_{(k)})=   \\
  \shoveleft {P\Big(\big(\sum_{u=1}^{k-1} L_{(u),1}\big)^{n_k-1}\prod_{u=1}^k
  \prod_{2\leq i<j\leq n_u} (L_{(u),j}-L_{(u),i}), \big((n_1-1,\ldots,1,0),\ldots} \\
\shoveright {\ldots,(n_{k-1}-1,\ldots,1,0),(n_k-2,\ldots,1,0)\big)\Big)=} \\
  \frac{(n_k-1)!}{\prod_{u=1}^{k-1} (n_u-1)!} P\Big(\prod_{u=1}^k \prod_{2\leq i<j\leq n_u}
                                                (L_{(u),j}-L_{(u),i}),\big((n_1-2,\ldots,1,0),\ldots,(n_k-2,\ldots,1,0)\big)\Big)=      \\
  (-1)^{\scriptscriptstyle \sum {n_u-1\choose 2}}\frac{(n_k-1)!}{\prod_{u=1}^{k-1} (n_u-1)!}
  \prod_{u=1}^{k} (n_u-1)!=(-1)^{\scriptscriptstyle \sum {n_u-1\choose 2}}(n_k-1)!^2.
\end{multline*}

 Hence, in the
expansion $-[D_\lambda]=c_1\sum L_{(1),i}+c_2\sum
L_{(2),i}+\ldots+c_k \sum L_{(k),i}$ we have
$$c_k=(n_k-1)!^2/\prod_{u=1}^k(n_u-1)!=(n_k-1)!/\prod_{u=1}^{k-1}(n_u-1)!$$
Proposition \ref{deg} then gives $\deg(D_\lambda)$ by substituting
(e.g.) $L_{(u),i}=0$ for $u<k$ and $L_{(k),i}=-1$ into
$[D_\lambda]$, which proves the Theorem.
\end{proof}

\section{ Some more explicit formulae for $GL(n)$.}
\label{andras}

Notice that the expression of $[D_\lambda]$ from \ref{short} is
homogeneous in $L_1,\ldots, L_n$ of degree 1. Moreover,
$[D_\lambda]\cdot \Delta$ is anti-symmetric of degree
$\deg\Delta+1$, hence it is the product of $\Delta$ and  a symmetric
polynomial of degree 1. Hence, by (\ref{degdlambda}), for some
constant $f_{\lambda}^{(n)}$ one has
\begin{equation}\label{flambda}
[D_\lambda]=-f_{\lambda}^{(n)}\cdot \sigma_1 \ \ \ \mbox{and} \ \ \
\ \deg(D_\lambda)=\frac{n}{|\lambda|}\cdot f_{\lambda}^{(n)},
\end{equation}
where $\sigma_1=\sum_{i=1}^nL_i$ as above. By similar argument one
deduces that for a free variable $t$ one has
$$ R_{\lambda}^{(n)}(t):= \sum_{\mu\in W\lambda}(\mu+t)\prod_{\beta\in
T_\mu}\frac{\mu+t+\beta}{-\beta}=A^{(n)}_\lambda
t+f^{(n)}_\lambda\sigma_1$$ for another constant $A^{(n)}_\lambda$.
Since the dual variety associated with the representation
$\rho_\lambda$ and the dual variety associated with the action
$\rho_\lambda\times \mbox{diag}$ of $GL(n)\times GL(1)$ are the
same, their degrees are the same too. Applying our main result for
these two groups and actions, one gets
\begin{equation}\label{tvar}
R^{(n)}_\lambda(t)=f^{(n)}_\lambda\cdot
\Big(\frac{n}{|\lambda|}\cdot t+ \sigma_1\Big).\end{equation} In
fact, (\ref{tvar}) is an entirely algebraic identity, and it is
valid for {\em any} weight $\lambda$ (i.e. not only for the dominant
weights).

Using (\ref{tvar}),  we prove some identities connecting different
weights.

\begin{lemma}\label{iesii}  For any $\lambda=
\sum_ia_iL_i$ define $\overline{\lambda}:=\sum_ia_iL_{n-i}$. Assume
that either (i) $\lambda'=\lambda+a\sigma_1$, or (ii)
$\overline{\lambda'}+\lambda=a\sigma_1$ for some $a\in \Z $. Then
$$\deg(D_{\lambda'})=\deg(D_{\lambda}).$$\end{lemma}

\begin{proof}
In case (i), using (\ref{tvar}), one gets
$T_{\lambda+a\sigma_1}=T_\lambda$ and
$f^{(n)}_{\lambda+a\sigma_1}=f^{(n)}_\lambda(\frac{na}{|\lambda|}+1)$,
hence
$f^{(n)}_{\lambda+a\sigma_1}/|\lambda+a\sigma_1|=f^{(n)}_\lambda/|\lambda|$.
Then apply (\ref{flambda}). For (ii) notice that $T_{\lambda'}=
\{-\overline{\beta}\,:\, \beta\in T_\lambda\}$, hence
$-[D_{\lambda'}]$ equals
$$\sum_{\mu\in W\lambda}
\overline{a\sigma_1-\lambda}\prod_{\beta\in
T_\mu}\frac{\overline{a\sigma_1-\lambda}+\overline{-\beta}}{\overline{\beta}}=
-\sum_{\mu\in W\lambda} (\lambda-a\sigma_1)  \prod_{\beta\in
T_\mu}\frac{\lambda-a\sigma_1+\beta}{-\beta}=$$
$$-f^{(n)}_{\lambda}\Big( -\frac{na}{|\lambda|}+1\Big)\sigma_1.$$
Therefore, $f^{(n)}_{\lambda'}=
f^{(n)}_\lambda(\frac{na}{|\lambda|}-1)$, or
$f^{(n)}_{\lambda'}/|\lambda'|=f^{(n)}_\lambda/|\lambda|$.

Geometrically the lemma  corresponds to the fact that tensoring with
the one dimensional representation or taking the dual representation
does not change the discriminant.
\end{proof}
It follows from the lemma that we  can always assume that the last
coefficient $a_n$ of $\lambda$ is zero.

\begin{example}\label{gammaab} {\bf The case of $\lambda=(a+b)L_1+b(L_2+\cdots + L_{n-1})$
with $a,b>0$ and $n\geq 3$.}\end{example}

\noindent  Set $\lambda':=(a+b)L_1+a(L_2+\cdots +L_{n-1})$. Then
$\overline{\lambda'}+\lambda=(a+b)\sigma_1$, hence
$\deg(D_\lambda)=\deg(D_{\lambda'})$ by \ref{iesii}. In particular,
$\deg(D_\lambda)$ is a {\em symmetric} polynomial in variables
$(a,b)$. In the sequel we deduce its explicit form. For this notice
that
$$T_\lambda=\{L_i-L_1,\ L_n-L_i: 2\leq
i\leq n-1\}\cup \{L_n-L_1\}.$$ It is convenient to write $\lambda$
as $aL_1-bL_n+b\sigma_1$. First, assume that $a\not=b$. Using
(\ref{flambda}) and Lemma~\ref{iesii} one gets
$$\deg(D_\lambda)=\frac{n}{|\lambda|}f^{(n)}_\lambda=\frac{n}{a+b(n-1)}
f^{(n)}_{aL_1-bL_n}\cdot
\Big(\frac{nb}{a-b}+1\Big)=\frac{n}{a-b}\cdot f^{(n)}_{aL_1-bL_n},$$
where $f^{(n)}_{aL_1-bL_n}\cdot \sigma_1$ equals
\begin{equation}\label{a-b}
\sum_{i\not=j}(aL_i-bL_j)\frac{aL_i-bL_j+L_j-L_i}{L_i-L_j}\prod_{k\not=i,j}
\frac{aL_i-bL_j+L_k-L_i}{L_i-L_k}\cdot\frac{aL_i-bL_j+L_j-L_k}{L_k-L_j}.\end{equation}
Clearly, $f^{(n)}_{aL_1-bL_n}=\lim_{L_n\to\infty}\,
(f^{(n)}_{aL_1-bL_n}\cdot \sigma_1)/L_n$. (In the sequel we write
simply $\lim$ for $\lim_{L_n\to \infty}$.) In order to determine
this limit, we separate $L_n$ in the expression $E$  of (\ref{a-b}).
We get three types of contributions: $E=I+II+III$, where $I$
contains those terms of the sum $\sum_{i\not=j}$ where $i,j\leq
n-1$; $II$ contains the terms with $i=n$; while $III$ those terms
with $j=n$.

One can see easily that $\lim  I$ is finite (for any fixed
$L_1,\ldots, L_{n-1}$), hence $\lim I/L_n=0$.

The second term (after re-grouping) is
$$II=\frac{S}{\prod\limits_{j\leq n-1}L_n-L_k},$$ where
$$S:=\sum_{j\leq n-1}\frac{\prod\limits_{k\leq
n-1}(aL_n-bL_j-L_n+L_k)(aL_n-bL_j+L_j-L_k)}
{\prod\limits_{n\not=k\not=j} L_k-L_j}.$$ By similar argument as
above, $S$ can be written as $S=\sum_{i=0}^n\ P_i\cdot L_n^{n-i}$,
where $P_i$ is a symmetric polynomial in $L_1,\ldots , L_{n-1}$ of
degree $i$ (and clearly also depends on $a,b$ and $n$). In
particular, $\lim II/L_n$ is the `constant'
$P_0=P_0(a,b,n)$.

Set $t_1:=aL_n-bL_j-L_n$ and $t_2:=aL_n-bL_j+L_j$. Notice that if we
modify the polynomial $P=\prod_{k\leq n-1}(t_1+L_k)(t_2-L_k)$ by any
polynomial situated in the ideal generated by the symmetric
polynomials $\sigma_1^{(n-1)}:=\sum_{i-1}^{n-1}L_i$,
$\sigma_2^{(n-1)}$, \ldots (in $L_1,\ldots , L_{n-1}$), then the
modification has no effect in the `constant' term $P_0$. In
particular, since $P=(t_1^{n-1}+\sigma_1^{(n-1)}t_1^{n-2}+\cdots)
(t_2^{n-1}-\sigma_1^{(n-1)}t_2^{n-2}+\cdots)$, in the expression of
$S$, the polynomial $P$ can be replaced by $t_1^{n-1}t_2^{n-1}$.
Therefore, if $[Q(t)]_{t^i}$ denotes the coefficient of $t^i$ in
$Q(t)$, then one has:
$$P_0(a,b,n)=\left[ \sum_{j\leq
n-1}\frac{(aL_n-bL_j-L_n)^{n-1}(aL_n-bL_j+L_j)^{n-1}}{\prod\limits_{n\not=k\not=j}
L_k-L_j}\right]_{L_n^n}.$$ Set
$$P_{a,b}(t)=At^2+Bt+C:=(at-t-b)(at-b+1).$$
Then
$$P_0(a,b,n)=\Big[P_{a,b}(t)^{n-1}\Big]_{t^n}\cdot
\sum_{j\leq n-1}\frac{L_j^{n-2}}{\prod\limits_{n\not=k\not=j}
L_k-L_j}.$$ Notice that the last sum is exactly $(-1)^{n-2}$ by
Lagrange interpolation formula.

A very similar computation provides $\lim III/L_n$, and using the
symmetry of $P_{a,b}$ one finds that
$$\deg(D_\lambda)=\frac{(-1)^nn}{a-b}\left(
\left[P_{a,b}(t)^{n-1}\right]_{t^n}-
\left[P_{a,b}(t)^{n-1}\right]_{t^{n-2}} \right)=$$
$$\frac{(-1)^nn}{a-b} \left( \left[P_{a,b}(t)^{n-1}\right]_{t^n}-
\left[P_{b,a}(t)^{n-1}\right]_{t^{n}} \right).$$ If $\partial$
denotes the `divided difference' operator $\partial
Q(a,b):=(Q(a,b)-Q(b,a))/(a-b)$, then the last expression reads as
$$\deg(D_\lambda)=(-1)^nn\cdot \partial \Big[P_{a,b}(t)^{n-1}\Big]_{t^n}.$$
By Corollary~\ref{fehpol}, this is valid for $a=b$ as well. By a
computation (using the multinomial formula):
\begin{equation}\label{positive}
\deg(D_\lambda)=n!(a+b-1)\sum_{i=1}^{[n/2]}\frac{A^{i-1}C^{i-1}(-B)^{n-2i}}
{i!(i-1)!(n-2i)!}.\end{equation} This has the following
factorization:
\begin{equation*}\label{newton}
\deg(D_\lambda)=n(n-1)(a+b-1)(-B)^{n-2[n/2]}\cdot
\prod_{i=1}^{[n/2]-1}(B^2+\xi_iAC),\end{equation*} where
$$\prod_{i=1}^{[n/2]-1}(t+\xi_i)=t^{[n/2]-1}+\frac{(n-2)!}{2!1!(n-4)!}t^{[n/2]-2}+
\frac{(n-2)!}{3!2!(n-6)!}t^{[n/2]-3}+\cdots.$$ E.g., for small
values of $n$ one has the following expressions for
$\deg(D_\lambda)$:
$$\begin{array}{ll}
n=3& 6(a+b-1)(-B)\\
n=4& 12(a+b-1)(B^2+AC)\\
n=5& 20(a+b-1)(-B)(B^2+3AC)\\
n=6& 30(a+b-1)(B^2+\xi AC)(B^2+\bar{\xi}AC), \ \mbox{where
$\xi^2+6\xi+2=0$}.\end{array}$$ For $n=3$, $\deg(D_\lambda)=
6(a+b-1)(2ab-a-b+1)$ is in fact the universal polynomial $F_{GL(3)}$
(which already was determined in Example~\ref{akbek} case $A_2$, and
from which one can determine the degrees of all the irreducible
$GL(3)$-representations by Corollary~\ref{fehpol}), compare also
with \cite{tevelev}[7.18].

 If we write
$u:=a-1$ and $v:=b-1$, then all the coefficients of the polynomials
$-B,A,C$, expressed in the new variables $(u,v)$, are non-negative.
Using (\ref{positive}), this remains true for $\deg(D_\lambda)$ as
well; a fact compatible with \cite{weyman}.

For $a=1$ the formula simplifies drastically (since $A=0$):
$\deg(D_\lambda)=n(n-1)b^{n-1}$. Symmetrically for $b=1$. For
$a=b=1$ we recover the well known formula $\deg(D_\lambda)=n(n-1)$
for the degree of the discriminant of the adjoint representation. (A
matrix is in the discriminant of the adjoint representation if it
has multiple eigenvalues i.e. the equation of the discriminant is
the discriminant of the characteristic polynomial.)

\begin{example}\label{abn} {\bf The case of $\lambda=aL_1+bL_2$ with
$a>b\geq 1 $ and $n\geq 3$.}\end{example}

\noindent     Since $T_\lambda=\{L_i-L_1,\ L_i-L_2\ :\ 3\leq i\leq
n\} \cup\{L_2-L_1\}$,  one has the following expression for
$f^{(n)}_\lambda\cdot \sigma_1$:
\begin{equation}\label{a+b}
\sum_{i\not=j}(aL_i+bL_j)\frac{aL_i+bL_j+L_j-L_i}{L_i-L_j}\prod_{k\not=i,j}
\frac{aL_i+bL_j+L_k-L_i}{L_i-L_k}\cdot\frac{aL_i+bL_j+L_k-L_j}{L_j-L_k}.
\end{equation}
We compute $f^{(n)}_\lambda$ as in \ref{gammaab}. For this we write
the expression $E$ from (\ref{a+b}) as $I+II+III$, where $I$, $II$
and $III$ are defined similarly as in \ref{gammaab}. The
computations of the first two contributions are similar as in
\ref{gammaab}: $\lim I/L_n=0$ while
$$\lim II/L_n=\Big[(at-t+b)^{n-1}(at+b-1)^{n-1}\Big]_{t^n}.$$
On the other hand, $III$ is slightly more complicated. For this
write $t_1:=aL_i+bL_n-L_i$, $t_2:=aL_i+bL_n-L_n$, and
$\Pi:=\prod_{n\not=k\not=i}L_i-L_k$. Moreover, we define the
relation $R_1\equiv R_2$ whenever $\lim R_1/L_n=\lim R_2/L_n$. Then
(for $t_1$ using the usual `trick' as above by separating terms from
the ideal ${\mathcal I}_{sym}$ generated by $\sigma_i^{(n-1)}$'s) we
have that $III=$
\begin{eqnarray*}
&=&\prod_{k\leq n-1}\frac{1}{L_n-L_k}\cdot \sum_{i\leq n-1}
\frac{\prod_{k\leq n-1}(t_1+L_k)\cdot (aL_i+bL_n+L_n-L_i)\cdot
\prod_{n\not=k\not=i}(t_2+L_k)}{\Pi} \\
&\equiv &\sum_{i\leq n-1} \frac{t_1^{n-1}}{\Pi\cdot L_n^{n-1}} \cdot
(aL_i+bL_n+L_n-L_i)\cdot
\prod_{n\not=k\not= i} (t_2+L_k)\\
&=&\sum_{i\leq n-1}\frac{t_1^{n-1}}{\Pi\cdot  L_n^{n-1}} \Big(
(b+1)L_n+(a-1)L_i\Big)
\frac{t_2^{n-1}+\sigma_1^{(n-1)}t_2^{n-2}+\cdots}{t_2+L_i}\\
&=& \sum_{i\leq n-1}\frac{t_1^{n-1}}{\Pi\cdot  L_n^{n-1}} \Big(
(b+1)L_n+(a-1)L_i\Big) \cdot
\end{eqnarray*}
$$\hspace{1cm}\cdot\Big(
\frac{t_2^{n-1}-(-L_i)^{n-1}}{t_2+L_i}
+\sigma_1^{(n-1)}\frac{t_2^{n-2}-(-L_i)^{n-2}}{t_2+L_i}+\cdots\Big)+$$
$$+ \sum_{i\leq n-1}\frac{t_1^{n-1}}{\Pi\cdot  L_n^{n-1}}
\Big( (b+1)L_n+(a-1)L_i\Big)
\frac{(-L_i)^{n-1}+\sigma_1^{(n-1)}(-L_i)^{n-2}+\cdots}{t_2+L_i}$$
From the first sum we can again eliminate terms from ${\mathcal
I}_{sym}$, while the second sum is a rational function in $L_n$,
where the numerator and denominator both have degree $n$, hence this
expression is $\equiv 0$. In particular,
$$III\equiv - \sum_{i\leq n-1}\frac{t_1^{n-1}}{\Pi\cdot L_n^{n-1}} \Big(
(b+1)L_n+(a-1)L_i\Big) \cdot
\frac{t_2^{n-1}-(-L_i)^{n-1}}{t_2+L_i}.$$
Therefore, $\deg(D_\lambda)$ equals
$$\frac{n}{a+b}\left[(at-t+b)^{n-1}(at+b-1)^{n-1}-
(bt+a-1)^{n-1}(bt+t+a-1)\cdot
\frac{(bt-t+a)^{n-1}-(-1)^{n-1}}{bt-t+a+1}\right]_{t^n}.
$$

For $n=3$ we recover the polynomial $F_{GL(3)}$ (cf. \ref{fehpol})
-- computed by \ref{gammaab} as well.

Also, for arbitrary $n$ but for $b=1$ the formula becomes simpler:
$$\deg(D_\lambda)=\frac{n}{(a+1)^2}\Big(
(n-1)a^{n+1}-(n+1)a^{n-1}+2(-1)^{n-1}\Big).$$ This expression
coincides with  Tevelev's formula from \cite{tevelevcikk}, see also
\cite[7.14]{tevelev}. (Notice that in \cite[7.2C]{tevelev}, Tevelev
provides a formula for arbitrary $a$ and $b$, which is rather
different from ours.)

\begin{example}\label{spec} {\bf Specializations of $\lambda=aL_1
+bL_2$.}\end{example}

\noindent Assume that $\lambda_s$ is a specialization of $\lambda$
(i.e. $\lambda_s$ is either $aL_1$ -- obtained by specialization
$b=0$, or it is $aL_1+aL_2$ -- by taking $b=a$). In this case, by
Corollary~\ref{fehpol} one has
$$
\deg(D_{\lambda_s})=\frac{\varepsilon(\lambda)|W_\lambda|}
{\varepsilon(\lambda_s)|W_{\lambda_s}|} \cdot (\,
\deg(D_\lambda)_{|\mbox{{\tiny specialized}}}).$$ In the first case,
$\deg(D_{aL_1+bL_2})|_{\,b=0}$ is $(-1)^nn(n-1)(a-1)^{n-1}$, and the
correction factor is $(-1)^n(n-1)$, hence we recover Boole's formula
$$\deg(D_{aL_1})=n(a-1)^{n-1}.$$
(Evidently, this can be easily deduced by a direct limit
computations -- as above -- as well.)

In the second case, $\deg(D_{aL_1+bL_2})|_{\,b=a}$ should be divided
by $-2$. We invite the reader to verify that this gives
$$\deg(D_{aL_1+aL_2})=\frac{n}{2a}\left[(t-1)(at+a-1)^{n-1}\cdot
\frac{(at-t+a)^{n-1}-(-1)^{n-1}}{at-t+a+1}\right]_{t^n}.$$ For $a=1$
we recover the well known formula of Holme \cite{holme} for the
degree of the dual variety of the Grassmannian $Gr_2(\C^n)$:
$$\deg(D_{L_1+L_2})=\frac{n}{2}\cdot \frac{1-(-1)^{n-1}}{2}.$$

\begin{example}
{\bf The case of $\lambda=L_1+L_2+L_3$, $n\geq 4$.} \end{example}

\noindent This case was studied by Lascoux \cite{lascoux}. Using
$K$-theory he gave an algorithm to calculate the degree of the dual
of the Grassmannian $Gr_3(\C^n)$ and calculated many examples. In
Proposition~\ref{gr3} we give a closed formula for the degree.

We write $\binom{n}{k}$ for the set of subsets of $\{1,2,\ldots,n\}$
with $k$ elements. For any $S\in\binom{n}{k}$ set $L_S:=\sum_{i\in
S}L_i$. Since $T_\lambda=\{L_i-L_j: \, i>3,\ j\leq 3\}$, one has
$$f^{(n)}_\lambda\cdot \sigma_1=\sum_{S\in\binom{n}{3}}
L_S\prod_{i\not\in S \atop j\in S} \, \frac{L_S+L_i-L_j}{L_j-L_i}.$$
Similarly as in the previous examples, we separate $L_n$, and we
write the above expression as the sum $I+II$ of two terms, $I$
corresponds to the subsets $S$ with $n\not\in S$, while $II$ to the
others. It is easy to see that $\lim_{L_n\to \infty}I/L_n=0$. In
order to analyze the second main contribution, it is convenient to
introduce the following expression (in variables $L_1,\ldots, L_n$
and a new free variable $t$):
$$R^{(n)}(t):=$$
$$\sum _{{\mathcal J}=\{j_1,j_2\}\in\binom{n}{2}}(t+L_{{\mathcal J}})
(t-L_{j_1})(t-L_{j_2})\prod_{i\not\in{\mathcal J}}\,\frac{
(t+L_{j_1}+L_i)(t+L_{j_2}+L_i)(L_{{\mathcal
J}}+L_i)}{(L_{j_1}-L_i)(L_{j_2}-L_i)}.$$ This is a homogeneous
expression in variables $(L,t)$ of degree $n+1$, hence it can be
written as
$$R^{(n)}(t)=\sum _{k=0}^{n+1}\, P_k^{(n)}t^{n+1-k},$$
where $P_k^{(n)}$ is a symmetric polynomial in variables $L$. The
point is that
$$II=\frac{R^{(n-1)}(L_n)}{\prod_{i\leq n-1}L_n-L_i}, \
\mbox{hence} \ \lim II/L_n=P_0^{(n-1)}.$$ Therefore,
$$f^{(n)}_\lambda=P_0^{(n-1)} \ \mbox{and} \
\deg(D_\lambda)=nP_0^{(n-1)}/3.$$ Next, we concentrate on the
leading coefficient $P^{(n)}_0$ of $R^{(n)}(t)$. For two polynomials
$R_1$ and $R_2$ of degree $n+1$, we write $R_1\equiv R_2$ if their
leading coefficients are the same. Set $t_r:=t+L_{j_r}$ for $r=1,2$;
and for each $k$ write (over the field $\C(L)$) $$
t_r^k=Q_{r,k}(t_r)(t_r+L_{j_1})(t_r+L_{j_2})+A_{r,k}t_r+B_{r,k}$$
for some polynomial $Q_{r,k}$ of degree $k-2$ and constants
$A_{r,k}$ and $B_{r,k}$. Then
$$\prod_{i\not\in{\mathcal J}}(t+L_{j_r}+L_i)=
\frac{t_r^n}{(t_r-L_{j_1})(t_r-L_{j_2})}+ \sum_{k=1}^n\sigma_k
Q_{r,n-k}+\sum_{k=1}^n\sigma_k\frac{A_{r,n-k}t_r+B_{r,n-k}}{
(t_r-L_{j_1})(t_r-L_{j_2})}.$$ In $P_0^{(n)}$ the first sum has no
contribution since it is in the ideal generated by the
(non-constant) symmetric polynomials, the second sum has no
contribution either, since its limit is zero when $t\to \infty$.
Therefore,  $R^{(n)}(t)$ is $\equiv $ with $$\sum _{{\mathcal
J}\in\binom{n}{2}}(t+L_{{\mathcal J}}) (t-L_{j_1})(t-L_{j_2})\cdot
\frac{(t+L_{j_1})^n}{(t+2L_{j_1})(t+L_{{\mathcal J}})}\cdot
\frac{(t+L_{j_2})^n}{(t+2L_{j_2})(t+L_{{\mathcal J}})}\cdot
\prod_{i\not\in{\mathcal J}}\,\frac{ L_{{\mathcal
J}}+L_i}{(L_{j_1}-L_i)(L_{j_2}-L_i)},$$ which equals
$$\sum _{{\mathcal  J}\in\binom{n}{2}}\frac{1}{t+L_{{\mathcal J}}}
\cdot \frac{(t-L_{j_1})(t+L_{j_1})^n}{t+2L_{j_1}}\cdot \frac{
(t-L_{j_2})(t+L_{j_2})^n}{t+2L_{j_2}}\cdot \prod_{i\not\in{\mathcal
J}}\,\frac{L_{{\mathcal J}}+L_i}{(L_{j_1}-L_i)(L_{j_2}-L_i)}.$$ Let
us define the---binomial-like---coefficients $\gbinom{n}{k}$ by the
expansion (near $t=\infty$):
\begin{equation}\label{expan}
\frac{(t-1)(t+1)^n}{t+2}=\sum_{k\leq n}\gbinom{n}{k}
t^k,\end{equation}
 Using
$$\frac{1}{t+1}=\frac{1}{t}-\frac{1}{t^2}+\frac{1}{t^3}-\cdots,$$
it provides the expansions
$$\frac{(t-L_{j_r})(t+L_{j_r})^n}{t+2L_{j_r}}=\sum_{k\leq n}\gbinom{n}{k}t^kL_{j_r}^{n-k}, \ \
\mbox{and}\ \ \frac{1}{t+L_{{\mathcal
J}}}=\frac{1}{t}-\frac{L_{{\mathcal J}}}{t^2}+\frac{L_{{\mathcal
J}}^2 }{t^3}-\cdots.$$ Therefore,
\begin{equation}\label{gr3formula}
P^{(n)}_0=\sum_{k\leq n, \, l\leq n\atop k+l\geq n+2}
(-1)^{k+l+n}\gbinom{n}{k}\gbinom{n}{l}{\mathcal
L}^{(n)}_{n-k,n-l},\end{equation}
 where
$${\mathcal L}^{(n)}_{n-k,n-l}:=
\frac{1}{2}\cdot \sum _{{\mathcal
J}=\{j_1,j_2\}\in\binom{n}{2}}L_{{\mathcal J}}^{k+l-n-2}
(L_{j_1}^{n-k}L_{j_2}^{n-l}+L_{j_1}^{n-l}L_{j_2}^{n-k})
 \cdot \prod_{i\not\in{\mathcal
J}}\,\frac{L_{{\mathcal J}}+L_i}{(L_{j_1}-L_i)(L_{j_2}-L_i)}.$$
(Here we already symmetrized ${\mathcal L}$ in order to be able to
apply in its computation the machinery of symmetric polynomials.
Using this index-notation for ${\mathcal L}^{(n)}_{n-k,n-l}$ has the
advantage that in this way this expression is independent of $n$, as
we will see later.) Next, we determine the constants $\gbinom{n}{k}$
and ${\mathcal L}^{(n)}_{n-k,n-l}$. The index-set of the sum of
(\ref{gr3formula}) says that we only need these constants for $k\geq
2$ and $l\geq 2$.

\vspace{1mm}

\noindent {\bf The constants $\gbinom{n}{k}$.} 
The identity
$$\frac{(t+1)^n(t-1)}{t+2}=(t+1)^n+3\sum_{i\geq 1}(-1)^i(t+1)^{n-i}$$
provides
\begin{equation}\label{cenk}
\gbinom{n}{k}=\binom{n}{k}+3\sum_{i\geq 1}(-1)^i\binom{n-i}{k}\ \
\mbox{for any $0\leq k\leq n$}.\end{equation} Notice that these
constants satisfy `Pascal's triangle rule':
$\gbinom{n}{k}=\gbinom{n-1}{k}+\gbinom{n-1}{k-1}$. This law together
with the `initial values' $\gbinom{n}{0}=1$ for $n$ even and $=-2$
for $n$ odd, and with $\gbinom{n}{n}=1$ determines completely all
$\gbinom{n}{k}$ for $0\leq k\leq n$. E.g., the first values are:

\vspace{2mm}

\begin{center}
\begin{tabular}{c||c|c|c|c|c|c|}
\hline  5& -2&-1&-2&1&2&1\\
\hline 4& 1& -2& 0& 1& 1&\\
\hline 3& -2& 0&0& 1& & \\
\hline 2& 1& -1& 1 &&& \\
\hline 1& -2& 1 &&&&\\
\hline 0& 1&&&&&\\
\hline \hline n/k & 0&1&2&3&4&5\\
\end{tabular}
\end{center}

\vspace{1mm}

\noindent {\bf The constants ${\mathcal L}^{(n)}_{n-k,n-l}$.} For
$a\geq 0$, $b\geq 0$, $n\geq 3$ and $a+b=n-2$ we set
$$X_{a,b}^{(n)}:=
\sum _{{\mathcal  J}=\{j_1,j_2\}\in\binom{n}{2}}
(L_{j_1}^{a}L_{j_2}^{b}+L_{j_1}^{b}L_{j_2}^{a})
 \cdot \prod_{i\not\in{\mathcal
J}}\,\frac{L_{{\mathcal J}}+L_i}{(L_{j_1}-L_i)(L_{j_2}-L_i)}.$$ Then
$X_{a,b}^{(n)}=X_{b,a}^{(n)}$. By the above homogeneity argument
$X_{a,b}^{(n)}$ is a constant. For $a>0$ and $b>0$ we separate $L_n$
and substitute $L_n=0$, and we get
\begin{equation}\label{xpascal}
X_{a,b}^{(n)}=X^{(n-1)}_{a-1,b}+X^{(n-1)}_{a,b-1}.\end{equation}
Therefore, $X_{a,b}^{(n)}$ can be determined from this `Pascal rule'
and the values $X_{n-1,0}^{(n)}$ (for $n\geq 3$). Next we compute
these numbers. We write $X_{n-2,0}^{(n)}$ as $I+II+III$, where $I$
corresponds  to ${\mathcal J}\in \binom{n-1}{2}$, and
$$II=\sum_{j\leq n-1}L_j^{n-2}\prod_{j\not=i\not=n}\frac{L_n+L_j+L_i}{(L_n-L_i)
(L_j-L_i)},$$
$$III=\sum_{j\leq
n-1}L_n^{n-2}\prod_{j\not=i\not=n}\frac{L_n+L_j+L_i}{(L_n-L_i)
(L_j-L_i)}.$$ Clearly, $X_{n-2,0}^{(n)}=\lim_{L_n\to 0}(I+II+III)$.
It is easy to see that $\lim I=0$ and $\lim II=1$. Moreover, using
 (\ref{tvar}), we get
$$III=\frac{L_n^{n-2}}{\prod_{i\leq n-1}(L_n-L_i)}\cdot \Big(
-\frac{1}{2}R_{2L_1}^{(n-1)}(L_n)+\frac{3}{2}L_n\sum_{j\leq
n-1}\prod_{ j\not=i\not=n}\frac{L_j+L_i}{L_j-L_i}\Big),$$ hence
\begin{equation}\label{initialx}
X^{(n)}_{n-2,0}=
1-\frac{n-1}{2}+\frac{3}{4}\Big(1+(-1)^n\Big).\end{equation} Now, we
return back to ${\mathcal L}^{(n)}_{n-k,n-l}$. The binomial formula
for $(L_{j_1}+L_{j_2})^{k+l-n-2}$ and (\ref{xpascal}) gives
$$2{\mathcal L}^{(n)}_{n-k,n-l}=\sum_{i=0}^{k+l-n-2}\binom{k+l-n-2}{i}
X^{(n)}_{l-2-i,n-l+i}=X^{(k+l-2)}_{l-2,k-2}.$$ Let
$\pmb\rangle^{(n)}_{a,b}$ be defined (for $a,b\geq 0$ and $a+b=n-2$)
by the Pascal rule and initial values $\pmb\rangle
^{(n)}_{n-2,0}=X^{(n)}_{n-2,0}$ and $\pmb\rangle^{(n)}_{0,n-2}=0$.
Symmetrically, define $\pmb\langle^{(n)}_{a,b}
=\pmb\rangle^{(n)}_{b,a}$, hence
$\pmb\rangle^{(n)}_{a,b}+\pmb\langle^{(n)}_{a,b} = X^{(n)}_{a,b}$.
It is really surprising that $X$ is another incarnation of the
constants $\gbinom{n}{k}$ (for $k\geq 1$). Indeed, comparing the
initial values of $\pmb\rangle^{(n)}_{a,b}$ and $\gbinom{n}{k}$ we
get that
$$\pmb\rangle^{(k+l-2)}_{l-2,k-2}=\gbinom{k+l-5}{k-1} \ \ (l\geq 2,\ k\geq 2).$$
Therefore, we proved the following fact.
\begin{proposition} \label{gr3} For any $n\geq 3$ and $2\leq k\leq n$ consider the
contants defined by (\ref{cenk}),
or by (\ref{expan}). Consider the weight $\lambda=L_1+L_2+L_3$  of
$GL(3)$ which provides the Pl\"ucker embedding of $Gr_3(\C^n)$. Then
the degree of the dual variety of $Gr_3(\C^n)$ is
%
$$\deg(D_\lambda)=\frac{n}{3}\cdot \mspace{-12mu} \sum_{\substack {\scriptscriptstyle
k\leq n-1,\\\scriptscriptstyle 4\leq l\leq n-1, \\
\scriptscriptstyle k+l\geq
n+1}}(-1)^{k+l+n-1}\gbinom{n-1}{k}\,\gbinom{n-1}{l} \, \gbinom{k
+l-5}{k-1} .$$

\end{proposition}

\bibliography{dual}
\bibliographystyle{amsplain}

\end{document}